\documentclass{article}
\usepackage[utf8]{inputenc}

\usepackage[english]{babel}
\usepackage[letterpaper,top=3cm,bottom=2cm,left=3cm,right=3cm,marginparwidth=1.75cm]{geometry}
\usepackage{amsmath}
\usepackage{amsfonts}
\usepackage{graphicx}
\usepackage[colorinlistoftodos]{todonotes}
\usepackage[colorlinks=true, allcolors=blue]{hyperref}

\title{A Computational Approach for Variational Integration of Attitude Dynamics on SO(3)}
\author{Nikhil Potu Surya Prakash \\
University of California, Berkeley \\
\textit{nikhilps@berkeley.edu}}
\date{}

\begin{document}

\maketitle
\begin{abstract}
In this article, a brief description of Discrete Mechanics and Variational Integrators which preserve the symplectic structure of the flow will be provided and a Newton-Raphson algorithm that can be used to solve implicit equations on the SO(3) manifold will be developed. These techniques will be used to simulate the rotational dynamics of a rigid body evolving on the Lie Group SO(3).   
\end{abstract}

\section{Introduction}
Discrete Mechanics first developed in [1] is a method of discretizing continuous differential equations for simulating on computers while preserving the structure of the manifold on which the dynamics evolve. The simplest discretization technique known is the Euler's technique which is used to numerically solve initial value problems by approximating the vector field to be a constant in a desired time interval 'h'. The Euler's method suffers from stability and other issues which have been to some extent addressed by a more accurate series of Runge-Kutta methods. Though the accuracy is increased, structures of the dynamics are not preserved in most cases. One such case is when a Runge-Kutta method is used to solve the Keppler problem, the energy and the angular momentum are not observed to be conserved. Discrete time models obtained via discrete mechanics are more desirable than other standard discretization schemes such as Euler’s step because they preserve certain invariance properties like kinetic energy, momentum, etc, of the system, and the computations can be done directly on the manifold, (because this discretization respects the manifold structure) thereby eliminating the problems associated with parametric representations. In the next sections, a brief description of Discrete Mechanics and its application to attitude dynamics of a rigid body on SO(3) will be provided. 
\section{Discrete Mechanics}
Consider a mechanical system with the configuration space $Q$
as a smooth manifold. Then the velocity vectors lie on the tangent bundle TQ of the manifold Q and the Lagrangian for the system can be defined as $L : TQ \rightarrow R$. In discrete mechanics, the velocity phase space TQ is replaced by $Q \times Q$ which is locally isomorphic to TQ. Let us consider an integral curve q(t) in the configuration space such that q(0) = q0 and q(h) = q1, where h represents the integration step. Then, the discrete Lagrangian $L_d : Q \times Q \rightarrow R$, which is an approximation of the action integral along the integral curve segment between q0 and q1, can be defined as
\begin{equation}\label{discrete lagrangian}
    L_d(q_0,q_1) \approx \int_0^h L(q(t),\dot q(t))dt
\end{equation}
Having defined the discrete Lagrangian, the action would now be the summation of all such Lagrangians along the path. The integral in the continuous case is replaced by the summation. This discrete Action is given by 
\begin{equation}\label{discrete action}
    A_d := \sum_{k=0}^{N-1} L_d(q_k,q_{k+1})
\end{equation}
where $q_i$ is the configuration of the system at $i$th time instant. \newline
Using the variation principles similar to the continuous case, the discretized equations of motion can be obtained as follows. [1] provides a detailed derivation of these equations.   
\begin{equation}\label{discrete EL}
     D_2L_d(q_{k-1},q_k)+D_1L_d(q_k,q_{k+1})=0 \; \forall \; k=0,1,...,N-1
\end{equation}
where $D_i$ is the derivative of the function with respect to the $i$th argument. \newline
Similarly, the discrete analogue of the Hamiltonian formulation can be obtained using the discrete Legendre transform. The continuous time Legendre transform is a map $\mathbb{F}L$ from the Lagrangian state space $TQ$ to the
Hamiltonian phase space $T^*Q$. Similarly, the discrete time Legendre transforms $\mathbb{F}^+L_d, \mathbb{F}^-L_d : Q \times Q \mapsto T^*Q$ [4] can be defined as
\begin{equation}\label{discrete EL1}
\mathbb{F}^{+}L_d(q_k,q_{k+1}) \mapsto (q_{k+1},p_{k+1}) = (q_k,D_2L_d(q_k,q_{k+1}))
\end{equation}
\begin{equation}\label{discrete EL2}
\mathbb{F}^{-}L_d(q_k,q_{k+1}) \mapsto (q_{k},p_{k}) = (q_k,-D_1L_d(q_k,q_{k+1}))
\end{equation}
where $p_i$ represents the corresponding conjugate momentum at $i$th time instant. \newline
In the next section, this routine will be used to develop the discretized equations of motion for the attitude dynamics of a rotating rigid body.

\section{Attitude dynamics using Discrete Mechanics}
In this section, a brief description of the attitude dynamics of a rigid body evolving on SO(3) manifold and its discrete equations of motion described in [3] will be presented. \newline
The kinetic energy of a rotating rigid body with angular velocity $\Omega$ is given by 
\begin{equation}\label{eq:KE}
    K = \frac{1}{2}\Omega^TJ\Omega = \frac{1}{2}tr(\widehat{\Omega}J_d\widehat{\Omega}^T) 
\end{equation}
where $\widehat{\Omega} \in \mathfrak{so(3)}$ is the skew symmetric 3$\times$3 tensor of $\Omega$ formulated according to \eqref{eq:vshom}, $J$ is the body moment of inertia matrix given by 
\begin{equation}
    J = \frac{1}{2}\int_\textit{B} \rho (X) \widehat{X}\widehat{X}^T d^3X 
\end{equation}
and $J_d$ is a matrix given by
\begin{equation}
    J_d = \frac{1}{2}\int_\textit{B} \rho (X) XX^T d^3X 
\end{equation}
The body moment of inertia $J$ is related to $J_d$ by
\begin{equation}\label{eq:JtoJd}
    J = tr(J_d)I_{3\times3}-J_d
\end{equation}
\eqref{eq:JtoJd} can be used to solve for $J_d$ and it can be obtained as
\begin{equation}
    J_d = \frac{1}{2} tr(J)I_{3\times3}-J
\end{equation}
The Lagrangian of the system would just contain the Kinetic Energy of the system and is given by
\begin{equation}
    L(R,\Omega) = K = \frac{1}{2}tr(\widehat{\Omega}J_d\widehat{\Omega}^T) 
\end{equation}
The rate of change of the rotation matrix $R$, which represents the orientation of the body fixed frame with respect to an inertial frame, is given by
\begin{equation}\label{eq:Rdot}
    \dot R = R \widehat{\Omega} 
\end{equation}
The Lagrangian of the system can be modified by replacing the angular velocity with $R^T \dot R$ from \eqref{eq:Rdot}
\begin{equation}\label{eq:AttLag}
    L(R,\dot R) = K = \frac{1}{2}tr(R^T \dot R J_d \dot R^T R)
\end{equation}
From the Lagrangian of the continuous time dynamics in \eqref{eq:AttLag}, the discrete Lagrangian can be obtained by approximating $\dot R$ using Euler's scheme with $h$ as the time step.
\begin{align}
    L_d(R_k,R_{k+1}) &\approx hL(R_k,\frac{R_{k+1}-R_k}{h}) \\
    &= \frac{h}{2}tr(\frac{R_k^T(R_{k+1}-R_k)}{h} J_d \frac{(R_{k+1}-R_k)^T R_k}{h}) \\
    &= \frac{1}{2h}tr((I_{3 \times 3}-F_k)J_d)
\end{align}
Here $R_i$ is the rotation matrix associated with the orientation of the body fixed frame with respect to the inertial frame at the $i$th time instant and $F_k = R_k^T R_{k+1}$ is used for notational brevity. Here it must be noted that $F_k$ is also a rotation matrix. \newline
The first order Hamilton's equations can be obtained using the above Lagrangian as follows. A detailed derivation of these equations can be found in [3].
\begin{equation}\label{eq:AttDyn}
Attitude \; Dynamics \;
\begin{cases}
\widehat{h\Pi} = F_kJ_d-J_dF_k^T \\
R_{k+1} = R_kF_k \\
\Pi_{k+1} = F_k^T\Pi_k+hu_k
\end{cases}
\end{equation}
where $\Pi_i$ is the conjugate momentum (the angular momentum in this case) and $u_i$ is the moment applied at the $i$th time instant.  \newline
The set of equations in \eqref{eq:AttDyn} can be used to simulate the system given the time series of the moments applied at each time instant, the initial orientation and the initial angular momentum. Since the propagation of the dynamics involves finding the rotation matrix $F_k$, the implicit matrix equation $\widehat{h\Pi} = F_kJ_d-J_dF_k^T$ needs to be solved at every time instant. A numerical technique will be developed in the next section to find the roots of the nonlinear matrix equation involving $F_k$.

\section{Solution to the Implicit Nonlinear Equation in $F_k$}
In this section, a Newton-Raphson type algorithm on the Lie Group SO(3) to solve the implicit nonlinear equation $\widehat{h\Pi} = F_kJ_d-J_dF_k^T$ to find $F_k \in SO(3)$ will be developed. Given the angular momentum $\Pi_k$ at the $k$th time instant, the implicit nonlinear matrix equation needs to be solved to find the incremental rotation matrix $F_k$ at every time instant to be used in the simulation of the dynamics.To preserve the manifold structure, it is essential that $F_k$ satisfies all the properties of a rotation matrix while solving the equation. Since $F_k$ is an incremental rotation matrix, there always exists a $\widehat{w} \in \mathfrak{so(3)}$, where $\mathfrak{so(3)}$ is the Lie Algebra of SO(3). A function $g(F_k)$ will be defined as follows, the zero of which will solve the required implicit nonlinear matrix equation.
\begin{equation}\label{eq:INME}
    g(F_k) = F_kJ_d-J_d F_k^T-\widehat{h\Pi_k}
\end{equation}
The above equation \eqref{eq:INME} can be expressed in terms of the vector in the lie algebra of $F_k$ as 

\begin{equation} \label{eq:INMELA}
    F(\widehat{w}) = e^{\widehat{w}}J_d-J_d (e^{\widehat{w}})^T-\widehat{h\Pi_k}
\end{equation}
It can be verified from the above equation \eqref{eq:INMELA} that the matrix $F(\widehat{w})$ is always a skew symmetric matrix. Therefore, a vector space homeomorphism can be established between the vector space of skew symmetric matrices and vector space $\mathbb{R}^3$ as follows

\begin{equation*}
\widehat{w} =
    \begin{bmatrix}
         w_1  \\
         w_2  \\
         w_3
    \end{bmatrix}^\times 
    = 
        \begin{bmatrix}
         0 & -w_3 & w_2   \\
         w_3 & 0 & -w_1  \\
         -w_2 & w_1 & 0
    \end{bmatrix}
\end{equation*}
and 
\begin{equation}\label{eq:vshom}
       \begin{bmatrix}
         0 & -w_3 & w_2   \\
         w_3 & 0 & -w_1  \\
         -w_2 & w_1 & 0
    \end{bmatrix}^\vee
    =
    \begin{bmatrix}
         w_1  \\
         w_2  \\
         w_3
    \end{bmatrix}
\end{equation}
Using this vector space homeomorphism, the skew symmetric matrix function of $\widehat{w}$ in \eqref{eq:INMELA} can be converted to a function of $w \in \mathbb{R}^3$ in $\mathbb{R}^3$.
\begin{equation}
    f(w) = [e^{\widehat{w}}J_d-J_d (e^{\widehat{w}})^T-\widehat{h\Pi_k}]^\vee
\end{equation}
The problem now boils down to finding a $w \in \mathbb{R}^3$ such that $f(w)=[0\;0\;0]^T$. A Newton-Raphson type algorithm can now be used to solve this using the jacobian of $f(w)$.But a closed form expression cannot be computed for neither $f(w)$ nor its jacobian $Df(w)$ as $f(w)$ is $F(\widehat{w})^\vee$ and hence the jacobians need to be computed jacobian of $F(\widehat{w})$. \newline
Since $F(\widehat{w})$ is always skew symmetric, its derivative with respect to $w_i$ is also skew symmetric. The derivative can be found to be
\begin{align}\label{eq:dFwdwi}
   \frac{\partial F(w)}{\partial w_i}  &= e^{\widehat{w}} \frac{\partial \widehat{w}}{\partial w_i} J_d-J_d (e^{\widehat{w}}\frac{\partial \widehat{w}}{\partial w_i})^T \\
    &= e^{\widehat{w}} \frac{\partial \widehat{w}}{\partial w_i} J_d+J_d \frac{\partial \widehat{w}}{\partial w_i}(e^{\widehat{w}})^T 
\end{align}
It can be verified that the derivative is skew symmetric and hence by using $\vee$ map defined above, the derivative $\frac{\partial f(w)}{\partial w_i}$ can be found as $[\frac{\partial F(w)}{\partial w_i}]^\vee$

\begin{equation}
   \frac{\partial f(w)}{\partial w_i}  = [\frac{\partial F(w)}{\partial w_i}]^\vee = [e^{\widehat{w}} \frac{\partial \widehat{w}}{\partial w_i} J_d+J_d \frac{\partial \widehat{w}}{\partial w_i}(e^{\widehat{w}})^T]^\vee
\end{equation}
These derivatives can be stacked to form the jacobian of $f(w)$ as
\begin{equation}
   Df(w) = [\frac{\partial f(w)}{\partial w_1} \; \frac{\partial f(w)}{\partial w_2} \; \frac{\partial f(w)}{\partial w_3}]^T  
\end{equation}
Though explicit expressions for $f(w)$ and $Df(w)$ were not obtained, it is still possible with above formulation to evaluate these at various $w$ and a numerical technique is feasible. \newline
With an initial guess $w^0$ and a step size of $\alpha$, the update equation in the Newton-Raphson algorithm is given by
\begin{equation}\label{eq:NRupdate}
   w^{n+1} = w^{n} - \alpha Df(w^n)^{-1}f(w^n)
\end{equation}
The vector $w^n$ is updated until $||f(w^n)||_2$ is satisfactorily close to 0 within user defined tolerance.

\section{Conclusion}
A brief description of Discrete Mechanics and Variational integrators was provided and its application to attitude dynamics of a rigid body evolving on the Lie Group SO(3) was provided. Mainly a Newton-Raphson type algorithm was developed to solve the implicit nonlinear matrix equation arising from Discrete Mechanics to obtain the incremental rotation matrix while preserving the manifold structure. 

\section*{MATLAB Package}
A package to simulate the attitude dynamics of a rigid body in MATLAB using the variational integration algorithm presented in this paper is available on request.

\end{document}